\documentclass[12pt]{article}
\usepackage[margin=1in]{geometry} 
\usepackage{setspace}
\usepackage[colorlinks]{hyperref}
\usepackage{mathtools}
\usepackage{caption}
\usepackage{yhmath,wasysym,skak}
\usepackage{amsmath, amsthm,amssymb,enumitem}
\usepackage{ytableau}
\usepackage{permute}
\usepackage{outlines}
\usepackage{tikz,marvosym}
\usepackage{pgfplots, txfonts}
\usepackage[marvosym]{tikzsymbols}
\usetikzlibrary{matrix,arrows.meta,shapes.geometric}
\usetikzlibrary{patterns,cd,fit,quotes}
\usetikzlibrary{arrows.meta}
\usetikzlibrary{intersections}
\usetikzlibrary{calc}
\usetikzlibrary{bending,decorations.markings}

\newtheorem{lemma}{Lemma}[section]

\newcommand{\DeclareAutoPairedDelimiter}[3]{%
  \expandafter\DeclarePairedDelimiter\csname Auto\string#1\endcsname{#2}{#3}%
  \begingroup\edef\x{\endgroup
    \noexpand\DeclareRobustCommand{\noexpand#1}{%
      \expandafter\noexpand\csname Auto\string#1\endcsname*}}%
  \x}
\DeclareAutoPairedDelimiter{\set}{\{}{\}}

\theoremstyle{definition}
\newtheorem{definition}[lemma]{Definition}
\newtheorem{example}[lemma]{Example}

\theoremstyle{remark}
\newtheorem{remark}[lemma]{Remark}

\theoremstyle{plain}
\newtheorem{prop}[lemma]{Proposition}
\newtheorem{theorem}[lemma]{Theorem}
\newtheorem{corollary}[lemma]{Corollary}

\def\la{\lambda}

\def\ol{\overline}
\def\x{\times}

\def\N{\mathbb{N}}
\def\C{\mathbb{C}}

\def\R{\mathbb{R}}

\def\mc#1{\mathcal{#1}}

\title{Lexicographic shellability of sects}

\begin{document}

\author{Aram Bingham\footnote{\texttt{aram@matmor.unam.mx}, Universidad de Chile}\hspace{5pt}   and N\'estor D\'iaz Morera\footnote{\texttt{ndiazmor@fitchburgstate.edu}, Fitchburg State University}}
\maketitle

% ABSTRACT
\begin{abstract}
     In this paper, we show that the Bruhat order on any sect of a symmetric variety of type $AIII$ is lexicographically shellable. Our proof proceeds from a description of these posets as rook placements in a partition shape which fits in a $p \times q$ rectangle. This allows us to extend an EL-labeling of the rook monoid given by Can to an arbitrary sect. As a special case, our result implies that the Bruhat order on matrix Schubert varieties is lexicographically shellable. 
\end{abstract}

\section{Introduction}

The Bruhat order on the symmetric group $\mc{S}_n$ is a fundamental object of study in algebraic combinatorics, containing important information about the representation theory of $GL_n(\C)$ as well as the geometry of its flag variety.  In this article, we consider the Bruhat order on a set of objects which are called \emph{$(p,q)$-clans} and which belong to an analogous  geometric and representation theoretic context. We find the notation and usage of \cite{yamamoto} the most amenable to our purposes.

\begin{definition}
Let $p$ and $q$ be two positive integers such that $p+q=n$. A \emph{$(p,q)$-clan} is an ordered set of $n$ symbols $c_1 \dots c_n$ such that: 
\begin{enumerate}
\item Each symbol $c_i$ is either $+$, $-$ or a natural number.
\item If $c_i \in \N$ then there is a unique index $j\neq i$ such that $c_i=c_j$.
\item The difference between the numbers of $+$ and $-$ symbols in the clan is equal to $p-q$; if $q>p$, then we have $q-p$ more minus signs than plus signs. 
\end{enumerate}
Clans are determined only up to equivalence based on the positions of matching pairs of natural numbers in the clan; the particular values of the natural numbers are not important. For example, ${+} 1 2 1 2 {-}$ and ${+} 2 3 2 3  {-}$ are equivalent $(3,3)$-clans. Also, $12{+}21$ is a $(3,2)$-clan and ${+}1{+}1$ is a $(3,1)$-clan. We denote the set of all $(p,q)$-clans by $\mc{C}_{p,q}$.
\end{definition}

By results of \cite{matsukiOshima}, the $(p,q)$-clans parameterize Borel subgroup orbits of a symmetric variety of type $AIII$. As symmetric varieties are \emph{spherical varieties} for which the Borel orbits are combinatorially parameterized, the set of Borel orbits also comes with a closure-containment order relation, or \emph{Bruhat order}, which can be combinatorially described. 

We will write $GL_n$ to indicate $GL(k^n)$ where $k$ is an algebraically closed field of characteristic not equal to 2. Fixing a basis for $k^{p+q}$, we denote by $B$ the Borel subgroup of $GL_{p+q}$ consisting of upper-triangular matrices. For clans $\gamma,\tau \in \mc{C}_{p,q}$, we have corresponding $B$-orbits $\mc{O}_\gamma$ and $\mc{O}_\tau$ in the symmetric variety $GL_{p+q}/(GL_p \x GL_q)$. Then the Bruhat poset $(\mc{C}_{p,q},\leq)$ is defined by the relation
\[\gamma \leq \tau \iff \mc{O}_\gamma\subseteq \ol{\mc{O}_\tau}\]
where $\ol{\mc{O}_\tau}$ denotes the Zariski closure of $\mc{O}_\tau$.

The Bruhat order on clans has been studied by Wyser in \cite{wyser16}, and the study of Bruhat order on symmetric varieties was initiated in greater generality in \cite{richardsonSpringer90}. Wyser raised the question of the shellability of this poset and noted that it does not satisfy other combinatorial properties common to Bruhat orders on flag varieties. In particular, $\mc{C}_{p,q}$ is neither thin nor Eulerian. However, shellability is a property with strong topological consequences, implying Cohen--Macaulayness and homotopy equivalence of the associated order complex to a wedge product of spheres. In the case of Bruhat orders, it has also been used to establish arithmetic Cohen--Macaulayness of certain Schubert varieties \cite{deconciniLakshmibai81}. Furthermore, lexicographic shellings of posets can be used to give combinatorial interpretation to the Möbius function (or Euler characteristic) of intervals or rank-selected subposets of shellable posets; see \cite{bjorner80, wachs} and references therein for background. 

Our main result, towards answering the shellability question for $\mathcal{C}_{p,q}$, says that slices of the poset of clans known as \emph{sects} and denoted $\mc{C}_{p,q}^\lambda$ are lexicographically shellable (and hence shellable) as they admit an EL-labeling; we review some of the background on poset shellability in Section~\ref{sec:shellability}.
\begin{theorem}\label{thm:main}
    The restriction of Bruhat order on $GL_{p+q}/(GL_{p} \x GL_{q})$ to any sect $\mc{C}_{p,q}^\lambda$ gives an EL-shellable poset. 
\end{theorem}
The sects can be described as preimages of Schubert cells under a natural projection map from $GL_{p+q}/(GL_p \x GL_q)$ to an associated Grassmannian variety $GL_{p+q}/P$. We remark that sect posets also fail to be thin and Eulerian in general; both failures can be observed within the sect poset $C_{2,2}^{(2,1)}$, which is isomorphic to a diamond with a largest element added above it. It is shown in \cite{binghamCan20, bingham21} that the largest sect of $\mc{C}_{p,q}$ is poset-isomorphic to the Bruhat order on matrix Schubert varieties, which can be defined as $B_p \x B_q$-orbit closures on the space of $p \x q$ matrices, where $B_p$ and $B_q$ are Borel subgroups of $GL_p$ and $GL_q$, respectively. The theorem above then implies the following, a special case of which is shown in \cite{can19}.

\begin{corollary}
    The Bruhat order on matrix Schubert varieties is EL-shellable.
\end{corollary}

Shellability of Bruhat posets of other type $A$ symmetric varieties ($SL_n/SO_n$ and $SL_{2n}/Sp_{2n}$, where orbits are parameterized by involutions and fixed-point-free involutions, respectively) was established in \cite{canTwelbeck, canCherniavskyTwelbeck}. The Bruhat poset $(\mc{C}_{p,q}, \leq)$ differs from these other type $A$ cases in that there is no longer a unique minimal element. The minimal elements of this poset are the \emph{matchless clans}: clans consisting only of ${+}$ and ${-}$ symbols. These correspond to the closed and minimum-dimensional orbits of the symmetric variety, which in turn correspond to Schubert cells of the related Grassmannian. As such, the full poset of $(p,q)$-clans is not bounded, though each matchless clan corresponds to a unique sect which is a bounded subposet $\mc{C}_{p,q}$, and the sects give a decomposition of the full poset. After encountering obstacles in our efforts to prove that arbitrary intervals of $\mc{C}_{p,q}$ are shellable, we turned to establishing the shellability of sects. This result is in fact independent from our original goal and generalizes several prior results, though we believe it also provides evidence of shellability of $\mc{C}_{p,q}$ and its maximal intervals.

Maximal intervals of $\mc{C}_{p,q}$ are close to being Bruhat (sub-)posets of involutions of $\mc{S}_n$. Shellability for these posets was established by Incitti in \cite{incitti04}. Our first attempts to demonstrate shellability of $\mc{C}_{p,q}$ made use of this fact, but we encountered several obstacles to translating Incitti's results to our context. Neither maximal intervals of $\mc{C}_{p,q}$ nor sects appear to arise naturally as intervals within known shellable posets, precluding use of the fact that shellability is inherited by intervals.  Further, we know of no geometric conditions guaranteeing shellability of a poset of stratifying subvarieties of a homogeneous variety, though we would find such conditions in the context of spherical varieties very desirable. Instead, our proof proceeds combinatorially based primarily on ideas from \cite{can19,bingham21}.

%%%%%%%%%%%%%%%%%%%%%%%%%%%%%%%%%%%%%%%%%%%%%%%%%%%
\section{Notation and preliminaries}

\subsection{Clans and rook placements}
The Bruhat order on $\mc{C}_{p,q}$ has been described in detail in \cite{wyser16, gandiniMaffei, bingham21}. Before we describe this order, we must precisely define the collections of clans known as ``sects."

Let $L=GL_{p}\x GL_q$ inside $G=GL_{p + q}$, so that $L$ can be realized as the subset of block-diagonal matrices.  This is the Levi factor of both the maximal parabolic subgroup $P$ of block upper-triangular matrices, as well as its opposite $P^-$. This gives projection maps 
\[G/P^- \longleftarrow G/L \longrightarrow G/P,\]
which are equivariant with respect to the $G$-action on each homogeneous space, and therefore equivariant with respect to the action of the Borel subgroup $B$ consisting of upper triangular matrices. This implies that the preimage of a $B$-orbit in $G/P$ (or $G/P^-$) is the union of a collection of $B$-orbits in $G/L$.  

The $B$-orbits of $G/P$ are called \emph{Schubert cells}, and they are known to be parameterized by the collection of partitions $\lambda$ that fit inside of a $p \x q$ rectangle; see for instance \cite{brion05}.  We denote this set of partitions by $\binom{[p+q]}{p}$, in reference to the fact that such a partition can be regarded as integer lattice paths inside the rectangle in the plane $\R^2$ with corners at the origin and the point $(q,p)$, encoded by a $p$-element subset of $[p+q]=\left\{1,2, \dots, p+q\right\}$.

\begin{definition}
    Let $C_\lambda$ denote the Schubert cell of $GL_{p+q}/P$  associated to the partition  $\lambda \in \binom{[p+q]}{p}$ where $P$ is as above. Then the \emph{sect} $\mc{C}_{p,q}^\lambda$ is the collection of clans $\gamma$ whose corresponding orbits satisfy $\pi(\mc{O}_\gamma)=C_\lambda$ where $\pi:G/L\to G/P$ is the natural projection map.
\end{definition}

The sects/Schubert cells/partitions in $\binom{[p+q]}{p}$ are in bijection with the matchless $(p,q)$-clans by constructing a lattice path from the origin to the point $(q,p)$ where $-$ signs correspond to ``east" steps and $+$ signs to ``north" steps. The partition $\lambda$ then consists of the boxes that lie inside the rectangle and weakly above this lattice path.

Let $\tau_\lambda$ denote the unique matchless $(p,q)$-clan associated to the partition $\lambda\in \binom{[p+q]}{p}$. The sect $\mc{C}_{p,q}^\lambda$ can then be described combinatorially as the set of clans $\gamma$ that become $\tau_\lambda$ under the following operation \cite{binghamCan20}. Label the symbols of $\gamma$ as $\gamma=c_1\cdots c_{p+q}$. For each pair $c_i=c_j\in \N$ with $i<j$, replace $c_i$ by a $-$ symbol and $c_j$ by a $+$ symbol. We refer to the resulting clan $\tau_\lambda$ as the \emph{base clan} of $\gamma$ and also denote it by $\tau_\gamma$ to indicate that it may be obtained from $\gamma$ via this procedure.

Our EL-labeling of $\mc{C}_{p,q}^\lambda$ is obtained by extending a labeling given in \cite{can19} (which itself is an extension of the labeling give by Edelman in \cite{edelman81}) and used to show that the so-called \emph{rook monoid} of partial permutations on $n$ letters is EL-shellable. We now introduce some rook placement notation and terminology.

\begin{definition}
    A \emph{rook placement of shape $\lambda$} is an assignment of rooks to the boxes of the partition of shape $\lambda$ in such a way that there is at most one rook in each row and each column. We denote the set of such rook placements by $R(\lambda)$. 
\end{definition}

\begin{prop}[{see \cite{bingham21}}] \label{prop:bij}
    The rook placements of shape $\lambda$ are in bijection with the clans of $\mc{C}_{p,q}^\lambda$ as follows. Let the positions of the $-$ symbols in $\tau_\gamma$ be ${i_1},\dots, {i_q}$ and the positions of the $+$ symbols be ${j_1}, \dots, {j_p}$ from left to right. Then for each pair $c_{i_k}=c_{j_l}\in \N$ in $\gamma$, we place a rook in the square with northeast corner $(k,l)$. 
\end{prop}

\begin{example}
 Let $\gamma$ be a $(3,3)$-clan given by $\gamma=1{+}{-}221$, with base clan $\tau_\gamma={-}{+}{-}{-}{+}{+}$. The labels on the positions of $\gamma$ are $c_{i_1}c_{j_1}c_{i_2}c_{i_3}c_{j_2}c_{j_3}$. Since $c_{i_1}=c_{j_3}$ and $c_{i_3}=c_{j_2}$, we have rooks in the boxes with northeast corners $(1,3)$ and $(3,2)$ as depicted in Figure~\ref{fig:one}. The presence of the gray rook below the partition shape and $\phi_\gamma$ will be explained in Definition~\ref{def:parper}.

 \begin{figure}[ht!]
    \centering
\begin{tikzpicture}[scale =.8]

\node(0) at (0,1){
\resizebox{.15\hsize}{!}{$
\gamma=1{+}{-}221
$}
};

\node (1) at (5,1) {
\resizebox{.1\hsize}{!}{$
\ytableausetup{mathmode, boxsize=1.2em}
\begin{ytableau}
 *(white) \symrook & *(white) &*(white)   \\
*(white)  & *(white) & *(white)\symrook    \\
*(white) &\none[ \textcolor{gray}{\symrook} ] 
 \\
\none[] &\none[\phi_{\gamma}=\textcolor{magenta!90!}{(3,1,2)}]
\end{ytableau}
$}
}; 

\draw [|->,line width=0.3mm] (0) edge node[below left,sloped] {} (1);
\end{tikzpicture}
  \caption{The bijection of Proposition~\ref{prop:bij} applied to a $(3,3)$-clan with base clan ${-}{+}{-}{-}{+}{+}$. Clans of this sect give rook placements of shape $\lambda=(3,3,1)$.}
    \label{fig:one}
\end{figure}
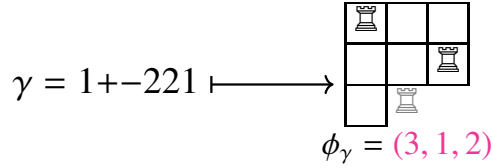
\end{example}

\begin{definition} 
    Write $[\lambda]$ to refer to the diagram of $\lambda$ as a collection of boxes (which can be regarded as a subset of $[q] \x [p]$). Given a rook placement $\rho \in R(\lambda)$, we can create the associated \emph{rank tableau} $rt_\rho: [\lambda] \to \N$ as a labeling of the boxes of $[\lambda]$ by the number of rooks weakly northwest of each box. 
\end{definition}
 
\begin{definition}
    We define a partial order $\preceq$ on $R(\lambda)$ by taking $\rho \preceq \pi$ if and only if 
    \[ rt_\rho(x) \leq rt_\pi(x)\]
    for every $x \in [\lambda]$. 
\end{definition}
\begin{theorem}[{\cite{bingham21}}] \label{thm:posetisom}
    The bijection of Proposition~\ref{prop:bij} defines a poset isomorphism of the Bruhat order on $(p,q)$-clans restricted to $\mc{C}_{p,q}^\lambda$ and $(R(\lambda), \preceq)$.
\end{theorem}
\begin{remark}
    In the special case $\lambda=(p^p)$, the theorem above gives a poset isomorphism of $\mc{C}_{p,q}^\lambda$ with the \emph{rook monoid} $R_p$, which is the Bruhat--Chevalley--Renner monoid of the algebraic monoid of square $p\x p$ matrices under matrix multiplication.
\end{remark}

Theorem~\ref{thm:posetisom} reduces the proof of Theorem~\ref{thm:main} to demonstrating the EL-shellability of $R(\lambda)$. However, to define our labeling, we must also refer to a \emph{partial permutation} $\phi_\gamma$ associated to a clan $\gamma$. To construct $\phi_\gamma$ we need a few auxiliary definitions. 

\begin{definition} We will have need to refer to certain patterns within clans.
\begin{itemize}
    \item A \emph{$1{+}{-}1$ pattern} in a clan $\gamma=c_1\cdots c_{p+q}$ is a collection of symbols $c_{i_1}c_{i_2}c_{i_3}c_{i_4}$, such that $i_1<i_2<i_3<i_4$, $c_{i_1}=c_{i_4}\in \N$, $c_{i_2}=+$ and $c_{i_3}=-$.
    \item A \emph{$1{2}1{2}$ pattern} in a clan $\gamma=c_1\cdots c_{p+q}$ is a collection of symbols $c_{i_1}c_{i_2}c_{i_3}c_{i_4}$, such that $i_1<i_2<i_3<i_4$, $c_{i_1}=c_{i_3}\in \N$ and $c_{i_2}=c_{i_4}\in \N$. 
    \item Analogously, define ${+}-$, ${+}11$, $11{+}$, $1122$, and $1221$ patterns.
\end{itemize}    
\end{definition}
\begin{definition}
    We will say that a $1{+}{-}1$ pattern $c_{i_1}c_{i_2}c_{i_3}c_{i_4}$ is \emph{simple} if there is no symbol $c_k$, $i_2<k<i_3$ where $c_k$ is a $+$ or $-$. 

    We will say that a $1{+}{-}1$ pattern with $c_{i_1}=c_{i_4}\in \N$ is \emph{innermost} if there is no $c_k=c_l\in \N$ with $i_1<k<i_2$ and $i_3<l<i_4$.
\end{definition}

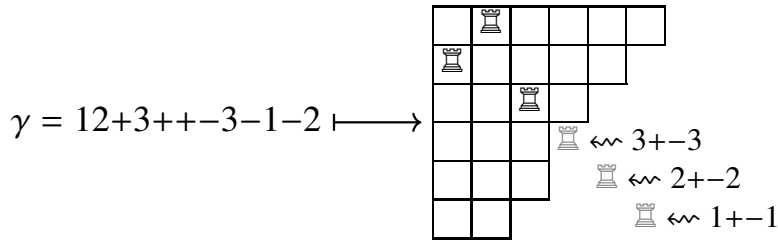
\begin{figure}[ht!]
    \centering
\begin{tikzpicture}[scale =.8]

\node(0) at (0,1){
\resizebox{.25\hsize}{!}{$
\gamma=12{+}3{+}{+}{-}3{-}1{-}2
$}
};

\node (1) at (7,1) {
\resizebox{.25\hsize}{!}{$
\ytableausetup{mathmode, boxsize=1.2em}
\begin{ytableau}
 *(white)  & *(white) \symrook  &*(white) & *(white) & *(white)& *(white)  \\
*(white)\symrook   & *(white) & *(white) &*(white) & *(white)  \\
*(white) & *(white) & *(white) \symrook & *(white) 
 \\
 *(white) & *(white) & *(white) &\none[ \textcolor{gray}{\symrook} ]  & \none[]{}&\none[\leftrightsquigarrow3{+}{-}3]
 \\
  *(white) & *(white) & *(white) &\none[]{} & \none[ \textcolor{gray}{\symrook}] & \none[]{}&\none[\leftrightsquigarrow 2{+}{-}2] 
 \\
  *(white) & *(white) & \none[]{} & \none[]{} & \none[]{} &\none[ \textcolor{gray}{\symrook} ] & \none[]{}&\none[\leftrightsquigarrow 1{+}{-}1]
 \\
%\none[] &\none[\phi_{\gamma}=\textcolor{magenta!90!}{(3,1,2)}]
\end{ytableau}
$}
}; 

\draw [|->,line width=0.3mm] (0) edge node[below left,sloped] {} (1);
\end{tikzpicture}
  \caption{Rook placement corresponding to a clan and its hidden rooks colored in gray. To the right of the hidden rooks, we indicate the symbols from the simple, innermost $1{+}{-}1$ pattern that gives rise to it. The hidden rook associated to the pattern $2{+}{-}2=\hat{c}_{i_2}\hat{c}_{j_2}\hat{c}_{i_5}\hat{c}_{j_5}$ is obtained from $\hat{\gamma}_1$, after changing the 1212 pattern to 1221 and deleting the symbols $3{+}{-}3=\hat{c}_{i_3}\hat{c}_{j_3}\hat{c}_{i_4}\hat{c}_{j_4}$. }
    \label{fig:hidden-rook}
\end{figure}

\begin{definition}\label{def:parper}
    The \emph{partial permutation} associated to a clan $\gamma=c_1\cdots c_{p+q}\in\mc{C}_{p,q}$ is the function $\phi_\gamma:[q]\to [p]\cup\{0\}$ defined algorithmically as follows. First label the positions of the $-$ symbols in $\tau_\gamma$ by ${i_1},\dots, {i_q}$ and the positions of the $+$ symbols as ${j_1}, \dots, {j_p}$ in ascending order. Then we read the symbols $c_1\cdots c_{p+q}$ left to right and construct $\phi_\gamma$ as follows.
    \begin{enumerate}
        \item If $c_{i_s}=c_{j_t}\in \N$ in $\gamma$, then $\phi_\gamma(s)=t$. These assignments represent rooks that are placed within the associated partition $\lambda$ under the bijection of Proposition~\ref{prop:bij}.
        \item After the previous step, we modify $\gamma$ by iteratively replacing all 1212 patterns by 1221 patterns to obtain a clan which we call $\hat{\gamma}_0 \in \mc{C}^\lambda_{p,q}$ and which has symbols $\hat{c}_{1}\cdots \hat{c}_{p+q}$.
        \item For each simple, innermost $1{+}{-}1$ pattern in $\hat{\gamma}_0$ of the form $\hat{c}_a \hat{c}_{j_l} \hat{c}_{i_k} \hat{c}_b$, set $\phi_\gamma(k)=l$. Then delete all of the symbols involved in any simple, innermost $1{+}{-}1$ pattern to obtain a new clan $\hat{\gamma}_1$ which inherits position labels from $\hat\gamma_0$.
        \item Repeat the procedure of the previous step on $\hat{\gamma}_1$, and so on until we obtain a clan $\hat{\gamma}_s$ which is free of $1{+}{-}1$ patterns.
        \item Once the clan $\hat\gamma_s$ which is free of $1{+}{-}1$ patterns is obtained, for any $k\in [q]$ which has not yet been assigned we let $\phi_\gamma(k)=0$.
    \end{enumerate}
    We will represent $\phi_\gamma$ using one-line notation. For instance, the partial permutation associated to the clan $\gamma=1{+}{-}221$ is $\phi_{\gamma}=(3,1,2)$ as shown in Figure~\ref{fig:one}.
\end{definition}
\begin{example}\label{ex:pperm}
    The clan $\gamma=12{+}{-}{+}{-}{-}{+}21\in\mc{C}_{5,5}$ has associated partial permutation $\phi_\gamma=(5,4,1,2,0)$. We note that for this clan, both $c_{i_2}c_{j_1}c_{i_3}c_{j_4}$ and $c_{i_2}c_{j_2}c_{i_4}c_{j_4}$ are simple, innermost $1{+}{-}1$ patterns. In step 3 of the above procedure, we would assign both $\phi_\gamma(3)=1$ and $\phi_{\gamma}(4)=2$ and then proceed to delete all six symbols involved in these patterns to obtain $\hat\gamma_1$, which is already free of $1{+}{-}1$ patterns. 

\end{example}
From this definition, we may obtain symbols in $\phi_\gamma$ that would correspond to rooks placed outside (below) the partition shape $\lambda$. We call these \emph{hidden rooks}. For instance, the full rook placement diagram for $\phi_\gamma=(5,6,4,3,2,1)$ from $\gamma=12{+}3{+}{+}{-}3{-}1{-}2$ is depicted in Figure~\ref{fig:hidden-rook}.

\begin{remark}\label{rk:pp}\hspace{1cm}\\

\begin{enumerate}[label=(\alph*)]
    \item Denote by $R_{q,p}$ the set of all partial injective transformations from $[q]$ to $[p]$; that is functions $\phi:[q] \to [p]\cup\{0\}$ whose restrictions to the set $N(\phi):=\left\{j\in[q] \mid \phi(j)\neq 0\right\}$ are injective. $R_{q,p}$ inherits a poset structure through bijection with the rook placements of $R((q^p))$ wherein we place a rook at the box with northeast corner $(k,l)$ for each $\phi(k)=l>0$. Since this poset is defined by northwest ``rank" conditions, the inclusion 
    \[(\mc{C}_{p,q}^\lambda\xrightarrow[\sim]{Proposition~\ref{prop:bij}} R(\lambda)) \hookrightarrow R_{(q^p)}\cong R_{q,p}\]
    is order-preserving. However, the image of this inclusion is not usually an interval of $R_{(q^p)}$, preventing us from applying the main result of \cite{can19} to $\mc{C}_{p,q}^\la$ to establish shellability. To observe this concretely, one may compare the partial permutations associated to the clans in Figure~\ref{fig:EL-sects} to their appearance in the full rook monoid $R_{3,3}$, as depicted in \cite[Figure 3.1]{can19}.

    \item The association of $\phi_\gamma\in R_{q,p}$ to $\gamma\in \mc{C}_{p,q}$ by Definition~\ref{def:parper} is one-to-one when restricted to a given sect $\mc{C}_{p,q}^\lambda$. One can see this by first noting that the association of Proposition~\ref{prop:bij} already injectively gives rook placements which can be regarded as partial permutations $\tilde{\phi}_\gamma\in R_{q,p}$ by the part (a) of this remark. The element $\phi_\gamma$ of Definition~\ref{def:parper} is then obtained from $\tilde{\phi}_\gamma$ by replacing some of the 0's in its one-line notation by other natural numbers.

    \item Finally, we mention that $R_{q,p}$ with the poset structure of part (a) is order isomorphic to the closure-containment (Bruhat) order on matrix Schubert varieties of the form $\ol{B_p^- w_{\phi} B_q}$, where $B_p^-$ denotes lower-triangular invertible $(p \x p)$-matrices, $B_q$ denotes upper-triangular invertible $(q\x q)$-matrices, and $w_\phi$ is the partial permutation matrix obtained by changing all the rooks in the associated rook placement to 1's, placing 0's elsewhere.
    \end{enumerate}
\end{remark}
Note that while Definition~\ref{def:parper} does not depend on the sect indexed by $\lambda$, the association will no longer be injective, and the labeling we will describe for our covering relations will fail to give an EL-shelling, if we do not restrict to a particular $\lambda$.

\subsection{Poset shellability} \label{sec:shellability}

The concept of lexicographic shellability of posets was pioneered by Bj\"orner \cite{bjorner80} and further developed in relation to Bruhat orders by Bj\"orner and Wachs \cite{bjornerwachs}. It continues to be widely studied due to its strong topological and algebraic consequences.

\begin{definition}
    Let $C(\mc{P}):=\{(u,v) \in \mc{P}\times \mc{P} \mid u \lessdot v \}$  denote the set of covering relations of a poset $(\mc{P}, \leq)$.  An \textbf{\emph{EL-labeling}} of $(\mc{P},\leq)$ is a map $\eta:C(\mc{P})\xrightarrow{} (\Lambda,\preceq)$, where $\Lambda$ is a totally ordered set, for which the following holds: 
    \begin{enumerate}
        \item[(1)] For each $u<v$, there is a unique (weakly) increasing sequence from $u$ to $v$. That is, there is a unique saturated chain $u\lessdot u_1 \lessdot \cdots \lessdot u_k \lessdot v$ with $\eta(u,u_1)\preceq \eta(u_1,u_2)\preceq \cdots \preceq \eta(u_k,v)$. 
        \item[(2)] The above  label sequence is lexicographically smaller than the label sequence for every
other saturated chain from $u$ to $v$.  That is, if $u\lessdot w<v$, with $w\neq u_1$ as defined above, then $\eta(u,u_1) \preceq \eta(u,w)$.
    \end{enumerate}
    A poset $\mc{P}$ is \textbf{\emph{EL-shellable}} or \textbf{\emph{lexicographically shellable}} if it admits an EL-labeling.
\end{definition}

\begin{example}
\begin{outline}
   \0 For permutations $u=u_1u_2\cdots u_n$ and $v=v_1v_2 \cdots v_n$ in the symmetric group $\mc{S}_n$ given in one-line notation, we say $v$ covers $u$ with respect to the Bruhat(--Chevalley) order, denoted by $u\lessdot_{\mathtt{BC}}v$, whenever $\ell(v)=\ell(u)+1$ where $\ell$ is the number of inversions, and 
        \1[(i)] $u_k=v_k$ for $k$ in $\{1,...,\widehat{i},...,\widehat{j},...,n\}$ (where the hat means we omit those numbers), 
        \1[(ii)]  $u_i=v_j$, $u_j=v_i$, and $u_i<u_j$.
    \0 If $C(\mc{S}_n)$ is the set of covering relations of $(\mc{S}_n,\leq_{\textbf{BC}})$ and $\Lambda=[n]\times [n]$ is the poset of pairs such that $(i,j) \leq (r,s)$ if $i<r$, or $i=r$ and $j<s$, we define $\eta(u,v):=(u_i,u_j)$. Then $\eta$ provides an EL-labeling for $\mc{S}_n$, as proposed by Edelman in \cite{edelman81}. In Figure~\ref{fig:EL-labeling}, we illustrate the corresponding EL-labeling for $\mc{S}_3$.
    \end{outline}
    
  \begin{figure}[ht]
    \centering
\begin{tikzpicture}[scale=0.7]   
\node (0) at (0,0) [above=0.2,draw] {$123$};
\node (1) at (-2,2) [above=0.2,draw] {$132$};
\node (2) at (2,2) [above=0.2,draw] {$213$};
\node (3) at (-2,6) [above=0.2,draw] {$231$};
\node (4) at (2,6) [above=0.2,draw] {$312$};
\node (5) at (0,8) [above=0.2,draw] {$321$};

\draw[red,line width=1mm] (0)--(2)--(3)--(5);

\path[-,line width=1.1] (0) edge node[below,sloped] {\small $(2,3)$} (1) edge node[below , sloped] { \small $(1,2)$} (2)
          (1) edge node[left] {\small $(1,2)$} (3)
              edge node[below right,sloped] {\small $(1,3)$} (4)
          (2) edge node[below left, sloped] {\small $(1,3)$} (3)
              edge node[right] {\small $(2,3)$} (4)
          (4) edge node[above, sloped] {\small $(1,2)$} (5)
          (3)  edge node[above,sloped] {\small $(2,3)$} (5);
\end{tikzpicture}
  \caption{EL-labeling of $\mc{S}_3$ with the unique increasing maximal chain highlighted in red.}
    \label{fig:EL-labeling}
\end{figure}
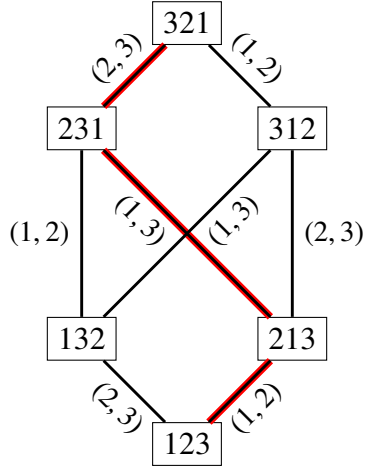
\end{example}

\begin{remark}
A finite, graded poset is \emph{shellable} if there is a well-ordering $<$ (called a \emph{shelling}) of the 	maximal chains of $\mc{P}$ such that if $\mathbf{c} <\mathbf{d}$, there exists another maximal chain $\mathbf{c}'<\mathbf{d}$  with $\mathbf{c} \cap \mathbf{d} \subseteq \mathbf{c}'\cap \mathbf{d}= \mathbf{d}\setminus \{x\}$ for some $x$.

A sufficient condition to ensure the shellability of a graded poset $\mc{P}$ is the EL-labeling property, as proved in \cite{bjorner80}. For more background on the theory of shellable posets and simplicial complexes, see \cite{wachs}.

\end{remark}

\subsection{Covering relations in sects}
We will precisely describe the covering relations of $\mc{C}_{p,q}^\lambda$. First, we introduce the notion of an underlying involution for a clan. Let the set of involutions in $\mc{S}_n$ be denoted by $\mc{I}_n$.

 \begin{definition}
\label{def:underlying}
	For a $(p, q)$-clan $\gamma = c_1 c_2 \cdots c_{n}$, the associated underlying involution $\pi_\gamma\in\mc{I}_{n}$ is defined as follows:
\begin{enumerate}
	\item $\pi_\gamma(i) = i$ if $c_i$ is either a $+$ or a $-$ for any $i$;
	\item  $\pi_\gamma(i) = j$ and $\pi(j)= i$ if $c_i = c_j \in \N$ is a matching pair of natural numbers for any $i, j$.
\end{enumerate}
%Note here that a different natural number is used for each pair of indices exchanged by the involution. 
 \end{definition} 
We can denote underlying involutions in notation similar to the clan notation by simply replacing all $+$ and $-$ signs by $\bullet$'s. For example, the underlying involution for the $(3,3)$-clan ${+}12{-}12$ would be $\bullet 12 {\bullet} 12$, or $156423$ in one-line notation.

 One of the important features of the description of the Bruhat order on $\mc{C}_{p,q}$ is that $\gamma\leq \tau$ implies that $\pi_\gamma \leq \pi_\tau$. That is, if clans are related in the Bruhat order on $\mc{C}_{p,q}$, then their underlying involutions are related in the Bruhat order on $\mc{I}_{p+q}$ as well. We now recall some notions from the Bruhat order on involutions that will be useful in the discussion that follows. 

\begin{definition}
    Let $\pi\in \mc{I}_n$, and consider $\pi:[n]\to[n]$ as a function, where $[n]=\left\{1,\dots,n\right\}$. A \emph{rise} of $\pi$ is a pair $(i,j)$ such that $i<j$ and $\pi(i)<\pi(j)$. For a clan $\gamma\in\mc{C}_{p,q}$, by a \emph{rise} of $\gamma$ we mean a rise of the underlying involution $\pi_\gamma$.
\end{definition}

\begin{definition}
    A rise $(i,j)$ of an involution $\pi$ is said to be \emph{free} if there is no $k$, $i<k<j$, such that $\pi(i)<\pi(k)<\pi(j)$. By a rise $(i,j)$ of a clan $\gamma$, we mean a rise of its underlying involution.
\end{definition}
The following is immediate from the definitions above.
\begin{lemma}
    A rise $(i,j)$ of a clan $\gamma=c_1\cdots c_{p+q}$ is free if and only if:
    \begin{itemize}
        \item There is no $k$ with $c_k=+$ or $c_k=-$ and $i<k<j$.
        \item There is no $c_k=c_l\in \N$ with $i<k<l<j$.
    \end{itemize}
    In other words, the only symbols between $c_i$ and $c_j$ are pairwise distinct natural numbers. 
\end{lemma}

We will show that the free condition guarantees that a rise is ``minimal" in the sense that free rises index covering relations of $\mc{C}_{p,q}$. First, we recall a few more notions from \cite{incitti04}.

\begin{definition}
    Let $i\in [n]$. For $\pi\in \mc{I}_n$, we will say that:
    \begin{itemize}
        \item $i$ is a \emph{fixed-point} of $\pi$ if $\pi(i)=i$;
        \item $i$ is an \emph{excedance} of $\pi$ if $\pi(i)>i$;
        \item $i$ is a \emph{deficiency} of $\pi$ if $\pi(i)<i$.
    \end{itemize}
\end{definition}
To each rise $(i,j)$ we can then associate a type $(a,b)$ where $a,b\in \left\{f,e,d\right\}$ according to whether $i$ and $j$ are fixed-points, excedances or deficiencies. For example, the involution $\bullet 12 \bullet 12=156423$ contains the rise $(1,2)$ which is of type $(f,e)$. Rises of type $(a,b)$ will also be called $ab$-rises. Incitti notes that for involutions, every $ee$-rise is associated with a symmetric $dd$-rise, $ef$-rises have corresponding $df$-rises, $fe$-rises correspond to $fd$-rises, and $ed$-rises correspond to $de$-rises. Thus, it suffices to consider only rises which are of type:
\[ (f,f),\ (f,e),\ (e,f),\ (e,e), \text{or } (e,d).\] 

\begin{definition}
    A rise $(i,j)$ is \emph{suitable} if it is free and of type $(f,f),\ (f,e),\ (e,f),\ (e,e)$, or $(e,d).$
\end{definition}

\begin{definition}
    For an involution $\pi\in\mc{I}_n$, we denote by $Inv(\pi)=\left\{(i,j)\mid i< j,\ \sigma(i)>\sigma(j)\right\}$ the set of inversions of $\pi$. We denote by $Exc(\pi)=\left\{i\in[n]\mid \sigma(i)>i\right\}$ the set of excedances of $\pi$. We denote the cardinalities of these sets by $inv(\pi)$ and $exc(\pi)$, respectively.
\end{definition}

The Bruhat poset of involutions is ranked (graded) by the \emph{length function} of \cite[Theorem 5.2]{incitti04},
\[\ell(\pi)= \frac{inv(\pi)+exc(\pi)}{2}.\]
On the other hand, the poset of clans is ranked using the length function of \cite[Definition 2.3.7]{yamamoto},
\begin{equation} \label{eq:length clans}
    \ell(\gamma)=\sum_{\substack{c_i=c_j\in \N \\ i<j }}j-i- \#\left\{a \in \N \mid c_s=c_t=a,\ s<i<t<j\right\}.
\end{equation}

\begin{lemma}
    The length of a clan is equal to the length of its underlying involution.
\end{lemma}
\begin{proof}
    First we note that the number of excedances of $\pi_\gamma$ is equal to the number of pairs of matching natural numbers that appear when we use the clan-like notation to represent $\pi_\gamma$. Next, observe that in every inversion $(i,j)$, at least one of $i$ or $j$ must be the position of a natural number. For each pair of natural numbers $c_i=c_j$, we can then define two sets of inversions,
    \[SInv(\pi;(i,j)) = \left\{(k,j) \in Inv(\pi) \mid k >i\right\},\]
    and
    \[FInv(\pi; (i,j))=\left\{(i,l) \in Inv(\pi) \mid \text{there is no $k$ with }(i,l)\in SInv(\pi; (k,l)\right\}.\]
    The condition on the inversions of $FInv(\pi;(i,j))$ ensures that these have not already been included among the inversions associated to another pair of mates. Thus, it is clear that 
    \[ inv(\pi)=\sum_{\substack{c_i=c_j\in \N \\ i<j }}FInv(\pi;(i,j))+SInv(\pi;(i,j)).\]
    To show that $\ell(\pi_\gamma)=\ell(\gamma)$, it then suffices to show that for a given pair of \emph{mates} $c_i=c_j\in \N$ in $\gamma$, we have that 
    \[\frac{FInv(\pi;(i,j))+SInv(\pi;(i,j))+1}{2}=j-i- \#\left\{a \in \N \mid c_s=c_t=a,\ s<i<t<j\right\}.\]
    That is, the contribution of a pair of mates to the length of the underlying involution is the same as the contribution of the mates in the clan.
    
    If we have a symbol $c_t$, $i<t<j$, which is a fixed point, then $(i,t)$ is an inversion. If $c_t$ is the second mate of another pair with first mate $c_s$, then it is again an inversion. However, it will not belong to $FInv(\pi;(i,j))$ if $s<i<t<j$. On the other hand, if $c_t$ is the first mate of a pair $c_t=c_u$ with $u<j$, then again $(i,t)$ is an inversion. The only other alternative is that $c_t=c_u$ is the first mate of a pair for $u>j$, in which case $(i,t)$ is not an inversion of $FInv(\pi;(i,j))$ but $(i,u)$ is. There can be no other inversions in $FInv(\pi;(i,j))$, so $FInv(\pi; (i,j))$ consists of exactly
    \[j-i- \#\left\{a \in \N \mid c_s=c_t=a,\ s<i<t<j\right\}\]
    inversions. 

    Next we consider $SInv(\pi;(i,j))$. For $i<k<j$, $(k,j)$ is an involution if $c_k$ is a fixed point, the first mate of a pair, or the second mate of a pair $c_l=c_k$ with $l>i$. There are $j-i-1$ such possible values for $k$, so we see that $SInv(\pi;(i,j))$ consists of \
    \[j-i-1-\#\left\{a \in \N \mid c_s=c_t=a,\ s<i<t<j\right\}\]
    elements. Adding to our count for $FInv(\pi; (i,j))$, this completes the proof.
\end{proof}

Now we can introduce the covering relations for $\mc{C}_{p,q}$. These are described implicitly in \cite[Theorem 2.8]{wyser16}, though we will refine this description when restricting this order to the sects. We will denote a covering relation of the Bruhat order on clans by $\gamma\lessdot \tau$.

\begin{lemma} \label{lem:covers}
    If $\gamma\lessdot \tau$ in the Bruhat order on $\mc{C}_{p,q}$, then $\tau$ is obtained from $\gamma$ by applying one of the following moves to a suitable rise of $\gamma$. The type of rise indicated refers to the type of the symbol indicated by the label beneath it.
    \begin{enumerate}
        \item ($ff$-rise) Replacing a pattern 
        $\underset{f}{+}\underset{f}{-} \mapsto 11, \text{ or } \underset{f}{-}\underset{f}{+}\mapsto 11.$
        \item ($fe$-rise) Replacing a pattern $\underset{f}{+}\underset{e}{1}1 \mapsto 1{+}1$ or $\underset{f}{-}\underset{e}{1}1\mapsto 1{-}1$.
        \item ($ef$-rise) Replacing a pattern $\underset{e}{1}1\underset{f}{+}\mapsto 1{+}1$ or $\underset{e}{1}1\underset{f}{-} \mapsto 1{-}1$.
        \item (non-crossing $ee$-rise) Replacing a pattern $\underset{e}{1}\underset{e}{2}12 \mapsto 1221$.
        \item (crossing $ee$-rise) Replacing a pattern $\underset{e}{1}1\underset{e}{2}2$ by either $1{+}{-}1$ (type $0$) or $1{-}{+}1$ (type $1$).
        \item ($ed$-rise) Replacing a pattern $\underset{e}{1}12\underset{d}{2}\mapsto 1212$.
    \end{enumerate}
\end{lemma}
\begin{proof}
    We know from \cite[Theorem 2.8]{wyser16} that to ascend in the Bruhat order in clans, one of the moves of the form listed above needs to be applied somewhere in $\gamma$. However, we must choose symbols such that the replacement only increases the length of the clan by one. Since the length of the clan agrees with the length of the underlying involution by the previous lemma, we can conclude from \cite[Theorem 5.1]{incitti04} that this move must occur at a suitable rise $(i,j)$, in which case the pattern replacement described above is exactly the ``covering transformation" $ct(i,j)$ described in \cite[Definition 3.2]{incitti04} on the level of underlying involutions. 
\end{proof}

\begin{lemma}\label{lem:sectcovers}
    If $\gamma\lessdot\tau$ in the Bruhat order on $\mc{C}_{p,q}^\lambda$, then $\tau$ is obtained from $\gamma$ by applying one of the following moves to a suitable rise of $\gamma$. 
       \begin{enumerate}
        \item ($ff$-rise) Replacing a pattern  $ \underset{f}{-}\underset{f}{+}\mapsto 11.$
        \item ($fe$-rise) Replacing a pattern  $\underset{f}{-}\underset{e}{1}1\mapsto 1{-}1$.
        \item ($ef$-rise) Replacing a pattern $\underset{e}{1}1\underset{f}{+}\mapsto 1{+}1$.
        \item (non-crossing $ee$-rise) Replacing a pattern $\underset{e}{1}\underset{e}{2}12 \mapsto 1221$.
        \item (crossing $ee$-rise) Replacing a pattern $\underset{e}{1}1\underset{e}{2}2$ by $1{+}{-}1$.
    \end{enumerate}
\end{lemma}
\begin{proof}
    By inspection, these are the moves from Lemma~\ref{lem:covers} that don't change the base clan $\tau_\gamma$. 
\end{proof}

%%%%%%%%%%%%%%%%%%%%%%%%%%%%%%%%%%%%%%%%%%%%%%%%%
\section{Shelling}
Permutation-like objects such as clans are amenable to shelling orders given by labeling edges in the poset diagram, since the elements can be represented by strings of constant length for which we ascend in the partial order by acting on the symbols of the string.  The strategy for shelling orders in this context is usually to provide a labeling on the edges so that the lexicographically smaller chains are those that act on the earliest possible symbols of the strings first.

One obstacle to finding an EL-labeling for the larger poset of clans is that the first symbols that can be acted upon in a given clan are sometimes not the first symbols of the clan. However, we can avoid this issue by labeling the covering relations within sects according to the change on the associated partial permutations. As we have mentioned in the introduction, this labeling extends the one given in \cite{can19} to rook placements of arbitrary shape. First we translate Lemma~\ref{lem:sectcovers} to the partial permutations of Definition~\ref{def:parper}. 

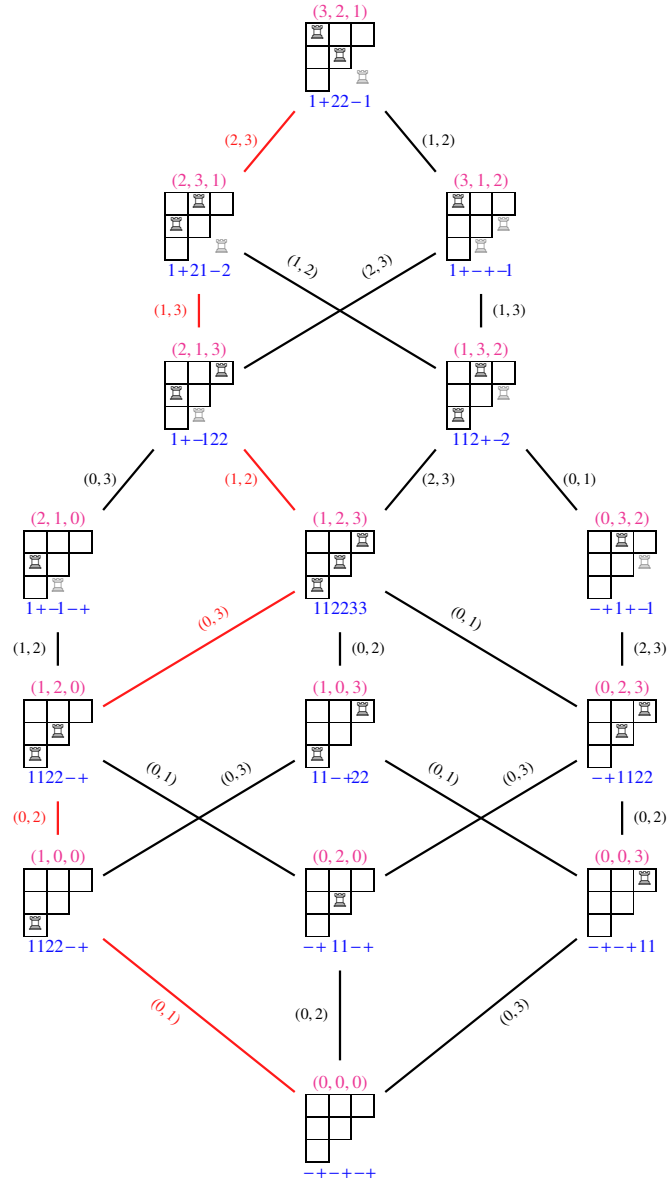
\begin{figure}[t!]
    \centering
  \begin{tikzpicture}[scale = 1.5]

%seventh  level 
\node (6) at (0,9+.5) {
\resizebox{.055\hsize}{!}{$
\ytableausetup{mathmode, boxsize=1.2em}
\begin{ytableau}
\none[] &\none[\textcolor{magenta!90!}{(3,2,1)}]
\\
 *(white)   \symrook & *(white) &*(white)   \\
*(white)  & *(white) \symrook \\
*(white) &\none[ ] &\none[ \textcolor{gray}{\symrook} ]
 \\
 \none[]& \none[\textcolor{blue!95!}{1\!+\!22\!-\!1}]
\end{ytableau}
$}
};

%sixth  level 
\node (52) at (1.25,7.5+.5) {
\resizebox{.055\hsize}{!}{$
\ytableausetup{mathmode, boxsize=1.2em}
\begin{ytableau}
\none[] &\none[\textcolor{magenta!90!}{(3,1,2)}]
\\
 *(white) \symrook & *(white) &*(white)   \\
*(white)  & *(white) &\none[ \textcolor{gray}{\symrook} ]   \\
*(white) &\none[ \textcolor{gray}{\symrook} ] 
 \\
 \none[]& \none[\textcolor{blue!95!}{1\!+\!-\!+\!-\!1}]
\end{ytableau}
$}
}; 
\node (51) at (-1.25,7.5+.5) {
\resizebox{.055\hsize}{!}{$
\ytableausetup{mathmode, boxsize=1.2em}
\begin{ytableau}
\none[] &\none[\textcolor{magenta!90!}{(2,3,1)}]
\\
 *(white)  & *(white) \symrook &*(white)   \\
*(white) \symrook & *(white) \\
*(white)  &\none[] & \none[ \textcolor{gray}{\symrook} ]  
 \\
 \none[]& \none[\textcolor{blue!95!}{1\!+\!21\!-\!2}]
\end{ytableau}
$}
};

%fifth level 
\node (42) at (1.25,6+.5) {
\resizebox{.055\hsize}{!}{$
\ytableausetup{mathmode, boxsize=1.2em}
\begin{ytableau}
\none[] &\none[\textcolor{magenta!90!}{(1,3,2)}]
\\
 *(white)  & *(white) \symrook &*(white)   \\
*(white)  & *(white) &\none[ \textcolor{gray}{\symrook} ]   \\
*(white) \symrook
 \\
 \none[]& \none[\textcolor{blue!95!}{112\!+\!-\!2}]
\end{ytableau}
$}
}; 
\node (41) at (-1.25,6+.5) {
\resizebox{.055\hsize}{!}{$
\ytableausetup{mathmode, boxsize=1.2em}
\begin{ytableau}
\none[] &\none[\textcolor{magenta!90!}{(2,1,3)}]
\\
 *(white)  & *(white)  &*(white) \symrook   \\
*(white) \symrook  & *(white)  \\
*(white)   &\none[ \textcolor{gray}{\symrook} ] 
 \\
 \none[]& \none[\textcolor{blue!95!}{1\!+\!-\!1 2 2}]
\end{ytableau}
$}
};

%fourth level 
\node (33) at (2.5,4.5+.5) {
\resizebox{.055\hsize}{!}{$
\ytableausetup{mathmode, boxsize=1.2em}
\begin{ytableau}
\none[] &\none[\textcolor{magenta!90!}{(0,3,2)}]
\\
 *(white)  & *(white) \symrook &*(white)   \\
*(white)  & *(white) &\none[ \textcolor{gray}{\symrook} ]   \\
*(white) 
 \\
 \none[]& \none[\textcolor{blue!95!}{-\!+\!1\!+\!-\!1}]
\end{ytableau}
$}
}; 
\node (32) at (0,4.5+.5) {
\resizebox{.055\hsize}{!}{$
\ytableausetup{mathmode, boxsize=1.2em}
\begin{ytableau}
\none[] &\none[\textcolor{magenta!90!}{(1,2,3)}]
\\
 *(white)  & *(white)  &*(white) \symrook   \\
*(white)  & *(white)\symrook  \\
*(white) \symrook
 \\
 \none[]& \none[\textcolor{blue!95!}{112233}]
\end{ytableau}
$}
}; 
\node (31) at (-2.5,4.5+.5) {
\resizebox{.055\hsize}{!}{$
\ytableausetup{mathmode, boxsize=1.2em}
\begin{ytableau}
\none[] &\none[\textcolor{magenta!90!}{(2,1,0)}]
\\
 *(white)  & *(white)  &*(white)   \\
*(white) \symrook   & *(white)  \\
*(white) &\none[ \textcolor{gray}{\symrook} ] 
 \\
 \none[]& \none[\textcolor{blue!95!}{1\!+\!-\!1\!-\!+}]
\end{ytableau}
$}
}; 

%third level 
\node (23) at (2.5,3+.5) {
\resizebox{.055\hsize}{!}{$
\ytableausetup{mathmode, boxsize=1.2em}
\begin{ytableau}
\none[] &\none[\textcolor{magenta!90!}{(0,2,3)}]
\\
 *(white)  & *(white)  &*(white) \symrook   \\
*(white)   & *(white)\symrook  \\
*(white) 
 \\
 \none[]& \none[\textcolor{blue!95!}{-\!+\!1122}]
\end{ytableau}
$}
};
\node (22) at (0,3+.5) {
\resizebox{.055\hsize}{!}{$
\ytableausetup{mathmode, boxsize=1.2em}
\begin{ytableau}
\none[] &\none[\textcolor{magenta!90!}{(1,0,3)}]
\\
 *(white)  & *(white)  &*(white)\symrook   \\
*(white)   & *(white)  \\
*(white) \symrook 
 \\
 \none[]& \none[\textcolor{blue!95!}{11\!-\!+\!22}]
\end{ytableau}
$}
};
\node (21) at (-2.5,3+.5) {
\resizebox{.055\hsize}{!}{$
\ytableausetup{mathmode, boxsize=1.2em}
\begin{ytableau}
\none[] &\none[\textcolor{magenta!90!}{(1,2,0)}]
\\
 *(white)  & *(white)  &*(white)   \\
*(white)   & *(white) \symrook  \\
*(white) \symrook 
 \\
 \none[]& \none[\textcolor{blue!95!}{1122\!-\!+}]
\end{ytableau}
$}
};

%second level 

\node (13) at (2.5,1.5+.5) {
\resizebox{.055\hsize}{!}{$
\ytableausetup{mathmode, boxsize=1.2em}
\begin{ytableau}
\none[] &\none[\textcolor{magenta!90!}{(0,0,3)}]
\\
 *(white)  & *(white)  &*(white) \symrook  \\
*(white)   & *(white) \\
*(white)  \\
 \none[]& \none[\textcolor{blue!95!}{-\!+\!-\!+\!11}]
\end{ytableau}
$}
};
\node (12) at (0,1.5+.5) {
\resizebox{.055\hsize}{!}{$
\ytableausetup{mathmode, boxsize=1.2em}
\begin{ytableau}
\none[] &\none[\textcolor{magenta!90!}{(0,2,0)}]
\\
 *(white)  & *(white)  &*(white)   \\
*(white)   & *(white) \symrook \\
*(white)  \\
 \none[]& \none[\textcolor{blue!95!}{\!-\!+11\!-\!+}]
\end{ytableau}
$}
};
\node (11) at (-2.5,1.5+.5) {
\resizebox{.055\hsize}{!}{$
\ytableausetup{mathmode, boxsize=1.2em}
\begin{ytableau}
\none[] &\none[\textcolor{magenta!90!}{(1,0,0)}]
\\
 *(white)  & *(white)  &*(white)   \\
*(white)   & *(white)  \\
*(white) \symrook \\
 \none[]& \none[\textcolor{blue!95!}{1122\!-\!+}]
\end{ytableau}
$}
};

%first level 
\node[] (0) at (0,0) {
\resizebox{.055\hsize}{!}{$
\ytableausetup{mathmode, boxsize=1.2em}
\begin{ytableau}
\none[] &\none[\textcolor{magenta!90!}{(0,0,0)}]
\\
 *(white)  & *(white)  &*(white)   \\
*(white)   & *(white)  \\
*(white) \\
 \none[]& \none[\textcolor{blue!95!}{-\!+\!-\!+\!-+}]
\end{ytableau}$}};

\draw [-,red!90!,line width=0.3mm] (0) edge node[below left,sloped] {\tiny$(0,1)$} (11);
\draw [-,line width=0.3mm] (0) edge node[left] {\tiny$(0,2)$} (12);
\draw [-,line width=0.3mm] (0) edge node[below right,sloped] {\tiny$(0,3)$} (13);

\draw [-,red!90!,line width=0.3mm] (11) edge node[left] {\tiny$(0,2)$} (21);
\draw [-,line width=0.3mm] (11) edge node[above right,sloped,xshift=10pt] {\tiny$(0,3)$} (22);
\draw [-,line width=0.3mm] (12) edge node[above left,sloped,xshift=-10pt] {\tiny$(0,1)$} (21);
\draw [-,line width=0.3mm] (12) edge node[above right,sloped,xshift=10pt] {\tiny$(0,3)$} (23);
\draw [-,line width=0.3mm] (13) edge node[above left,sloped,xshift=-10pt] {\tiny$(0,1)$} (22);
\draw [-,line width=0.3mm] (13) edge node[right] {\tiny$(0,2)$} (23);

\draw [-,line width=0.3mm] (21) edge node[left] {\tiny$(1,2)$} (31);
\draw [-,red!90!,line width=0.3mm] (21) edge node[above right,sloped] {\tiny$(0,3)$} (32);
\draw [-,line width=0.3mm] (22) edge node[right] {\tiny$(0,2)$} (32);
\draw [-,line width=0.3mm] (23) edge node[above left,sloped] {\tiny$(0,1)$} (32);
\draw [-,line width=0.3mm] (23) edge node[right] {\tiny$(2,3)$} (33);

\draw [-,line width=0.3mm] (31) edge node[left] {\tiny$(0,3)$} (41);
\draw [-,red!90!,line width=0.3mm] (32) edge node[left] {\tiny$(1,2)$} (41);
\draw [-,line width=0.3mm] (32) edge node[right] {\tiny$(2,3)$} (42);
\draw [-,line width=0.3mm] (33) edge node[right] {\tiny$(0,1)$} (42);

\draw [-,red!90!,line width=0.3mm] (41) edge node[left] {\tiny$(1,3)$} (51);
\draw [-,line width=0.3mm] (41) edge node[above right,sloped,xshift=10pt] {\tiny$(2,3)$} (52);
\draw [-,line width=0.3mm] (42) edge node[above left,sloped,xshift=-10pt] {\tiny$(1,2)$} (51);
\draw [-,line width=0.3mm] (42) edge node[right] {\tiny$(1,3)$} (52);

\draw [-,red!90!,line width=0.3mm] (51) edge node[left] {\tiny$(2,3)$} (6);
\draw [-,line width=0.3mm] (52) edge node[right] {\tiny$(1,2)$} (6);

\end{tikzpicture}

    \caption{EL-shellability of $\mc{C}_{3,3}^{(3,2,1)}$ with the unique maximal increasing chain highlighted in red.}
    \label{fig:EL-sects}
\end{figure}

\begin{lemma}\label{lem:ppcovers}
    Let $\gamma=c_1\cdots c_{p+q}\in\mc{C}_{p,q}^\lambda$, with symbols labeled as in Proposition~\ref{prop:bij} and with partial permutation $\phi_\gamma=(a_1,\dots,a_q)$. Then:
    \begin{enumerate}
        \item An $ff$-rise move applied to $(c_{i_k},c_{j_l})$ changes the symbol $a_k$ from $0$ to $l$.
        \item An $fe$-rise move applied to $(c_{i_k},c_{i_m})$ where $c_{i_m}=c_{j_l}\in \N$ swaps the symbol pair $(a_k,a_m)$ from $(0,l)$ to $(l,0)$.
        \item An $ef$-rise move applied to $(c_{i_k},c_{j_m})$ where $c_{i_k}=c_{j_l}\in \N$ changes the symbol $a_k$ from $l$ to $m$.
        \item A non-crossing $ee$-rise move applied to $(c_{i_k},c_{i_m})$ where $c_{i_k}=c_{j_l}\in \N$ and $c_{i_m}=c_{j_n}\in \N$ swaps the symbol pair $(a_k,a_m)$ from $(l,n)$ to $(n,l)$.
        \item A crossing $ee$-rise move applied to $(c_{i_k},c_{i_m})$ where $c_{i_k}=c_{j_l}\in \N$ and $c_{i_m}=c_{j_n}\in \N$ swaps the symbol pair $(a_k,a_m)$ from $(l,n)$ to $(n,l)$.
     \end{enumerate}
\end{lemma}
\begin{proof}
    This is straightforward to check by inspecting the partial permutation associated to the clan which results from each type of move identified in Lemma~\ref{lem:sectcovers}.
\end{proof}
We now define our labeling.
\begin{definition}\label{def:stand.label}
    Suppose that $\gamma\lessdot\tau$ is a covering relation where $\gamma$ and the covering moves are as in the notation of Lemma~\ref{lem:ppcovers}. Then the \textbf{\emph{standard labeling}} $L$ is the one that labels the covering relation $(\gamma,\tau)\in C(\mc{C}_{p,q}^\lambda)$ with an element of $\N^2$ (totally ordered lexicographically) as follows.
    \begin{enumerate}
        \item If $\tau$ is obtained from $\gamma$ by an $ff$-rise move, then apply the label $(0,l)$.
        \item If $\tau$ is obtained from $\gamma$ by an $fe$-rise move, then apply the label $(0,l)$.
        \item If $\tau$ is obtained from $\gamma$ by an $ef$-rise move, then apply the label $(l,m)$.
        \item If $\tau$ is obtained from $\gamma$ by a non-crossing $ee$-rise move, then apply the label $(l,n)$.
         \item If $\tau$ is obtained from $\gamma$ by a crossing $ee$-rise move, then apply the label $(l,n)$.
    \end{enumerate}
\end{definition}
An example of this labeling applied to one of the sects of $\mc{C}_{3,3}$ is depicted in Figure~\ref{fig:EL-sects}. The label applied to the covering relation $1{+}21{-}2\lessdot 1{+}22{-}1$ is $(2,3)$, since our partial permutation changes from $(2,3,1)$ to $(3,2,1)$ via a non-crossing $ee$-rise move.

We next proceed to identify the lexicographically smallest chain in any interval in a given $\mc{C}_{p,q}^\lambda$. Suppose $\gamma <\tau$ in $\mc{C}_{p,q}^\lambda$, and let their associated partial permutations be denoted $\phi_\gamma=(a_1,\dots, a_q)$ and $\phi_\tau=(b_1,\dots, b_q)$. We will also denote 
\begin{align*}
    D(\gamma,\tau)&=\left\{j \in [q]\mid a_j\neq b_j\right\},\\
    D_0(\gamma,\tau)&=\left\{j \in [q]\mid a_j=0\text{ and }b_j\neq 0 \right\}.
\end{align*}
\begin{definition}
    Define the \emph{action index} $t$ for $\gamma<\tau$ as follows.
    \begin{itemize}
        \item If $D_0$ is empty, then $t$ is the index such that $a_t=\min\left\{a_i \mid i \in D(\gamma,\tau)\right\}$ ($a_t$ is the ``lowest" out-of-place rook).
        \item If $D_0$ is not empty, then it contains a least element, $k_0$. First we construct the auxiliary set $E_0$ of \emph{entry points} as follows: we say $k$ is an entry point if $k\geq k_0$, $a_k=0$, there are more non-zero symbols among ${b_1,\dots,b_k}$ than there are among ${a_1,\dots,a_k}$, and there is a free $fe$-rise $(i_k,i_{m})$ or there is a free $ff$-rise $({i_k},j_r)$.\footnote{The entry points are the analog of the set $P(x,y)$ from \cite[p. 273]{can19}.} For any entry point $k$, there is a clan that covers $\gamma$ (and is below $\tau$) obtained by applying one of these rise moves. Note that the entry point $i_k$ may be part of both a free $ff$-rise and a free $fe$-rise. 
        \begin{enumerate}
            \item If the number of symbols $a_j$ in $\phi_\gamma$ that are non-zero equals the number of symbols $b_j$ in $\phi_\tau$ that are non-zero, then our action index is the value $k$ indexing the entry point associated to the suitable $fe$-rise $(i_k,i_m)$ where $c_{i_m}=c_{j_v}$, and the associated covering relation has the smallest possible label $(0,v)$.
            \item If the number of symbols $a_j$ in $\phi_\gamma$ that are non-zero is less than the number of symbols $b_j$ in $\phi_\tau$ that are non-zero, then our action index is the $k$ indexing the entry point associated to the suitable $ff$-rise $(i_k,j_v)$ or $fe$-rise $(i_k, i_m)$  with the smallest label $(0,v)$.
        \end{enumerate}
    \end{itemize}
\end{definition}

\begin{example}\label{ex:act.ind}
    Let $\gamma=12{+}{-}{+}{-}{-}{+}21$ and let $\tau=12{+}3{+}{-}{-}321$ (so $\tau$ is the maximal clan of the sect of $\gamma$). From Example~\ref{ex:pperm} we have $\phi_\gamma=(5,4,1,2,0)$, and one may check that $\phi_\tau=(5,4,3,2,1)$. From this, we obtain $D(\gamma,\tau)=\{3,5\}$ and $D_0(\gamma,\tau)=\{5\}$. Furthermore, there are fewer non-zero symbols in $\phi_\gamma$ than in $\phi_\tau$, so the action index is $t=k_0=5$, associated to the $ff$-rise $(i_{5},j_3)$. 

    If $\gamma'=12{+}{-}{+}{-}3321$ with $\phi_{\gamma'}=(5,4,1,2,3)$, then $D_0(\gamma',\tau)$ is empty, so in this case $t=3$ since $\phi_{\gamma'}$ and $\phi_\tau$ first differ at the third position.
\end{example}
\begin{definition}
    Define the \emph{cover value} $v$ for $\gamma<\tau$ as follows, where the action index $t$ has $a_t=x$. 
    \begin{itemize}
     \item If $D_0$ is empty, then we let $v$ be the minimum possible value of the set $\left\{x+1,\dots,p\right\}$ such that either:
    \begin{enumerate}[label=(\alph*)]
    \item there is no $u\in[q]$ with $a_u=v$, or
    \item there is an index $u>t$ with $a_u=v$, and some index $s\in D(\gamma,\tau)$ such that $s\geq u$ and $b_s= a_t <a_s$. We note that the existence of this index $s$ is necessary to guarantee that after applying the covering move with label $(a_t,v)$, the resulting clan is still beneath $\tau$.
\end{enumerate}
        \item If $D_0$ is not empty, then the action index is defined as the position of the entry point that allows us to chose the smallest possible covering label $(0,v)$. This $v$ is then our cover value.
    \end{itemize}
\end{definition}
\begin{example}\label{ex:cov.val}
    Continuing with the clans of Example~\ref{ex:act.ind}, the cover value for the pair $\gamma, \tau$ is $v=3$ since $D_0(\gamma,\tau)$ is non-empty and the $ff$-rise associated to the action index is of the form $(i_t,j_v)=(i_5,j_3)$. On the other hand, the pair $\gamma',\tau$ has cover value $v=2$ since (following part (b) of the case where $D_0$ is empty) the index $u=4$ satisfies $a_4=2$ with $s=5$ providing $b_5=1=a_3<a_5=3$. 
\end{example}

\begin{lemma} \label{lem:key}
    If $\gamma <\tau$ is an arbitrary order relation of $\mc{C}^\lambda_{p,q}$, then there is a unique clan $mcc_\tau(\gamma)$ called the \textbf{minimal covering clan} for the interval $[\lambda,\tau]$, denoted $mcc_\tau(\gamma)$ and satisfying:
    \begin{enumerate}
        \item $\gamma\lessdot mcc_\tau(\gamma)\leq \tau$
        \item $L(\gamma,mcc_\tau(\gamma))=(a_t,v)$ where $t$ and $v$ are the action index and cover value defined above, and $L(\cdot,\cdot)$ denotes the standard labeling of Definition~\ref{def:stand.label}.
    \end{enumerate}
\end{lemma}
\begin{proof} The definitions of $t$ and $v$ already determine a partial permutation corresponding to $mcc_\tau(\gamma)$; it remains to verify that there is a valid clan covering move that achieves this partial permutation. We handle this with cases depending on $D_0$.

\noindent\textbf{($D_0$ empty)}  If $a_t=x$ and $v=x+k$ for some $k\geq 1$, then by definition of the cover value $v$, either there is $u>t$ with $a_u=v$ or there is no symbol with this value in $\phi_\gamma$. We consider these separately. 
    \begin{enumerate}
             \item If $a_u=v$ for some $u>t$, then one of the following is true by construction of the partial permutations.
             \begin{enumerate}
                 \item $c_{i_t}c_{i_u}c_{j_x}c_{j_{v}}$ is a $1212$ pattern to which a non-crossing $ee$-rise move can be applied.
                 \item $c_{i_t}c_{j_x}c_{i_u}c_{j_{v}}$ is a $1122$ pattern to which a crossing $ee$-rise move can be applied.
                 \item $c_{i_t}c_{j_x}c_{j_{v}}c_{i_u}$ is a $11{+}{-}$ pattern, and there is a simple $1{+}{-}1$ pattern $\hat{c}_a \hat{c}_{j_{x+1}}\hat{c}_{i_u}\hat{c}_d$ in $\hat{\gamma}_0$. In this case, an $ef$-rise move can be applied to the suitable rise $(i_t,j_{v})$.
                 \item $c_{j_x}c_{i_t}c_{i_uc_{j_{v}}}$ is a ${+}{-}11$ pattern, and there is a simple $1{+}{-}1$ pattern $\hat{c}_a\hat{c}_{j_x}\hat{c}_{i_t}\hat{c}_d$ in $\hat{\gamma}_0$. In this case, an $fe$-rise move can be applied to the suitable rise $(i_t,i_u)$.
                  \item $c_{j_x}c_{i_t}c_{j_{v}}c_{i_u}$ is a ${+}{-}{+}{-}$ pattern, and there are simple $1{+}{-}1$ patterns $\hat{c}_{a_1} \hat{c}_{j_x} \hat{c}_{i_t}\hat{c}_{d_1}$ and $\hat{c}_{a_2} \hat{c}_{j_{v}} \hat{c}_{i_u}\hat{c}_{d_2}$ in $\hat{\gamma}_0$. In this case, an $ff$-rise move can be applied to the suitable rise $(i_t,j_{v})$.
             \end{enumerate}
            \item If there is no $u\in[q]$ with $a_u=v$, then $c_{i_t}c_{j_x}c_{j_{v}}$ is a $11{+}$ pattern for which $(i_t,j_{v})$ is a suitable $ef$-rise, and there is no $1{+}{-}1$ pattern in $\hat{\gamma}_0$ containing $\hat{c}_{j_{x+1}}$.         
            \end{enumerate}
             In each case, the appropriate rise move defines a clan $mcc_\tau(\gamma)\gtrdot\gamma$, and $L(\gamma,mcc_\tau(\gamma))=(a_t,v)=(a_t,a_{t}+k)$ by construction. 
             That $mcc_\tau(\gamma)\leq \tau$ follows from \cite[Lemma 4.2]{can19}, since the map from clans to partial permutations is order-preserving (see Remark~\ref{rk:pp}, part (a)) and we have $\phi_{mcc_\tau(\gamma)}\leq \phi_\tau$. Essentially, one may restrict to the non-zero part of the partial permutation, and the Deodhar description of Bruhat order on partial permutations ensures that after interchanging $a_t$ and $a_{t}+k$, we remain below $\tau$ which has $b_t>a_t$.
             
   \noindent \textbf{($D_0$ non-empty)} In this case, our action index $t$ gives us $a_t=0$. We let $mcc_\tau(\gamma)$ be the clan obtained by applying the appropriate covering move to the free $ff$-rise $(i_t, j_v)$ or $fe$-rise $(i_t,i_m)$. The label given to this covering relation will be $(0,v)$, and the fact that $mcc_\tau(\gamma)\leq\tau$ follows from \cite[Lemma 4.3]{can19}, once again comparing on the level of partial permutations. 
  \end{proof}

The following is immediate from the construction of $mcc_\tau(\gamma)$, based on the definitions of $t$ and $v$.
\begin{lemma} \label{lem:atom labels}
Under the standard edge-labeling, $L(\gamma, mcc_\tau(\gamma))$ has the smallest label among the labels on edges to atoms of the interval $[\gamma,\tau]$.    
\end{lemma}
Theorem~\ref{thm:main} follows immediately from our next result, which completes the argument that our standard labeling provides a lexicographic shelling of $\mc{C}_{p,q}^\lambda$.
\begin{theorem}
    The standard edge-labeling is an EL-labeling of $\mc{C}^\lambda_{p,q}$. 
\end{theorem}
\begin{proof}
   Let $[\gamma,\tau]$ be an interval in $\mc{C}_{p,q}^\lambda$ for a given valid $\lambda$. Then consider the maximal chain $\gamma=\gamma_0\lessdot \gamma_1\lessdot \dots \lessdot \gamma_n=\tau$ where $\gamma_{i+1}=mcc_\tau(\gamma_i)$. By Lemma~\ref{lem:atom labels}, this chain is lexicographically smallest among the maximal chains in $[\gamma,\tau]$. Furthermore, its labels are weakly increasing since the minimal covering clan is defined at each stage by:
   \begin{enumerate}
       \item First moving the smallest non-zero symbols of $\phi_{\gamma_i}$ into the non-zero positions of $\phi_\tau$, or introducing non-zero symbols of $\phi_{\gamma_i}$ as necessary, giving priority to entry points where this can be done with the smallest possible label.
       \item Then swapping non-zero symbols to ascend in Bruhat order,  beginning with the smallest out-of-place symbols.
   \end{enumerate}
   Next we argue that any other maximal chain will have a descent among its labels. If $\gamma\lessdot \rho\neq \gamma_1$ is another covering relation, then $L(\gamma,\rho)=(l_1,l_2)>L(\gamma,\gamma_1)=(a_t,v)$. We can argue that in every chain of $[\rho, \tau]$, there is a covering relation with label that is smaller than $(l_1,l_2)$, so that any maximal chain $\gamma<\rho<\dots <\tau$ must not be weakly increasing. We will divide this into cases according to whether $l_1>a_t$, or $l_1=a_t$ and $l_2>v$. 

   \begin{enumerate}
       \item \textbf{($l_1>a_t$)} If $a_t\neq 0$, then we know that $a_t<b_t$, so at some point this symbol will need to be increased to reach $b_t$ in the partial permutation of $\tau$. Any covering move that increases $a_t$ has first label equal to $a_t$, so in any saturated chain of the interval $[\rho,\tau]$ there is a covering relation with first label $a_t<l_1$, creating a descent. In case $a_t=0<l_1$, then there is a position $k_0\leq t$ with $a_{k_0}=0\neq b_{k_0}$. In any saturated chain of $[\rho,\tau]$ there is then a covering relation with first label $0<l_1$ required to increase the value of $a_{k_0}$, creating a descent.

       \item \textbf{($l_1=a_t$, $l_2>v$)} We again need subcases.
       \begin{enumerate}
           \item First suppose that $a_t\neq 0$, implying that $D_0$ is empty. Then if $u$ is the index such that $a_u=v$ and $y$ is the index such that $a_y=l_2$, it must be the case that $t<y<u$ since otherwise the rise corresponding to the covering relation $\gamma\lessdot \rho$ would not be free. 
       
        Denote $\phi_\rho=(r_1,\dots, r_q)$. Now $r_y=a_t$ and $r_u=a_u=v$. By definition of $v$, there is an index $s\geq u>y$ such that $b_s=a_t <a_u=v$. If $s=u$, then in order to move the rook corresponding to $a_u=v$ down to height $a_t$, we need to apply a covering move with label $(a_t, v)<(l_1,l_2)$ at some point in any saturated chain of $[\rho,\tau]$, and we are done.       
        
        If, on the other hand, $s\neq u$, then by the minimality of $v$ among the set $\left\{a_t+1,\dots q\right\}$ we must have $a_s> a_u$. The symbol $a_s=r_s$ is not affected by the covering move from $\gamma$ to $\rho$, which implies that in each saturated chain of the interval $[\rho,\tau]$, there must be some covering move with first label $b_s= a_t=l_1$ in order to move the corresponding rook into place. But now $r_u=a_u<a_s=r_s$ implies that there is no free rise move which can swap a rook at height $a_t=r_y$ with a rook at height $r_s$ directly, without first interchanging rooks at heights $r_y=a_t$ and $r_u=v$. Thus, there must be a covering relation with label $(a_t,v)<(a_t,l_2)=(l_1,l_2)$ at some point in every saturated chain of $[\rho,\tau]$, finishing this case as well. 

        \item Now suppose that $a_t=0$. The only covering moves with label $(0,l)$ are $ff$-rises and $fe$-rises. We consider these subcases separately. 
        
        \begin{enumerate}
            \item Suppose first that $\gamma\lessdot \gamma_1$ comes from an $fe$-rise move applied to the rise $(i_t,i_u)$ with $c_{i_u}=c_{j_v}$. 
        Denote $\phi_{\gamma_1}=(d_1,\dots,d_q)$ and $\phi_\rho=(r_1,\dots, r_q)$ as before. 
        
        If $r_t\neq 0$, then the only possibility is that $\rho$ is obtained from a different $fe$-rise move applied to a rise $(i_t, i_m)$ with $c_{i_m}=c_{j_{l_2}}$ for $m<u$ and $l_2>v$. Since $\gamma_1<\tau$, there must be at least as many non-zero symbols among the first $m$ symbols of $\tau$ as there are in $\gamma_1$. The number of non-zero symbols among $\left\{r_1,\dots, r_{m}\right\}$ is one less than the number of non-zero symbols among $\left\{d_1,\dots,d_m\right\}$, implying that there must be a covering move with first label $0$ which brings a rook into column $m$ in any saturated chain of $[\rho,\gamma]$. For $\rho$, $(0,v)$ remains the minimal possible label among its entry points and is still only obtainable through an $fe$-rise move. Since $(i_m,i_u)$ is now a free rise for $\rho$, the only way to bring a rook into column $m$ is via another $fe$-rise move with label $(0,z)$ for some $v\leq z<l_2$, creating a descent.

        If $r_t=0$, then $\rho$ is obtained by applying a rise move at a different entry point $k>t$ with label $(0,l_2)$. Since $\phi_\rho$ needs at least one more rook than it has among the first $t$ columns of its diagram in order to reach $\tau$, $t$ must be used as an entry point for at least one covering relation in every saturated chain of $[\rho,\tau]$. If the label on this covering move is $(0,z)$ for $z<l_2$, then we are done. Otherwise, it is $(0,w)$ for some $fe$-rise $(i_t, i_m)$ with $t<m<u$ and $c_{i_m}=c_{j_w}$, and we call the clan reached after applying this move $\sigma$. Then we can apply reasoning as in the previous paragraph before to conclude that $(i_m,i_u)$ will become a free $fe$-rise for the interval $[\sigma,\tau]$ leading to a descent of the form $(0,z)$ with $v\leq z<l_2<w$.

        \item Finally, we consider the case where $\gamma\lessdot \gamma_1$ comes from an $ff$-rise move. Let this rise be $(i_t,j_v)$.

        If $r_t\neq0$, then once again $\rho$ is obtained from an $fe$-rise move applied to the rise $(i_t, i_m)$, now with $c_{i_m}=c_{j_{l_2}}$ for $i_m<j_v<j_{l_2}$. Arguing similarly to the previous case, at some point in any saturated chain of $[\rho, \tau]$, a rook must be introduced in column $m$. The only way to achieve this is now via the $ff$-rise move corresponding to free rise $(i_m,j_v)$, or via another $fe$-rise move for a rise $(i_m, i_u)$ with $c_{i_u}=c_{j_z}$ and $z<{l_2}$. In either case, we obtain a descent.

        If $r_t=0$, the argument is similar to the situation where $\gamma_1$ is obtained from $\gamma$ by an $fe$-rise. The only difference is that in case $t$ is used as an entry point for an $fe$-rise $(i_t,i_m)$ with label $(0,w)$ and $w>l_2$ in order to reach clan $\sigma$, then either $(i_m,j_v)$ is a free $ff$-rise or there is a free $fe$-rise $(i_m,i_u)$ with $c_{i_u}=c_{j_z}$, $z<l_2$, one of which must occur in the interval $[\sigma,\tau]$ (as in the previous paragraph). 
    
        \end{enumerate}
        \end{enumerate}

   \end{enumerate}
   This completes the argument that $\gamma\lessdot \gamma_1 \lessdot \dots \lessdot \gamma_n= \tau$ is the unique weakly increasing saturated chain among $[\gamma,\tau]$, and hence the proof that $L$ is an EL-labeling of this interval. Since the standard labeling gives an EL-shelling of any interval of $\mc{C}_{p,q}^\lambda$, this implies the poset is EL-shellable.
\end{proof}

\subsection*{Acknowledgements}

We thank the anonymous referee for several helpful suggestions which improved the clarity of our exposition. We are also grateful to Mahir Can for introducing us to this research area and posing the kinds of questions that led to this article. 

%BIBLIOGRAPHY

\end{document}